\documentclass[journal]{IEEEtran}
\usepackage{amsfonts}
\usepackage[comma,numbers,square,sort&compress]{natbib}

\newcommand{\be}{\begin{equation}}
\newcommand{\ee}{\end{equation}}
\newcommand{\bee}{\begin{equation*}}
\newcommand{\eee}{\end{equation*}}
\newcommand{\bea}{\begin{eqnarray}}
\newcommand{\eea}{\end{eqnarray}}

\usepackage{graphicx}          

\usepackage[dvips]{epsfig}    
\usepackage{latexsym,amssymb}
\usepackage{amsmath, amsbsy}
\usepackage{amsopn, amstext}

\newtheorem{remark}{Remark}
\newtheorem{theorem}{Theorem}
\newtheorem{definition}{Definition}
\newtheorem{lemma}{Lemma}

\newtheorem{proposition}{Proposition}

\usepackage{CJK}
\begin{document}

\title{Integrated Stabilization Policy over a Software Defined Network}



\author{Cheng Tan \IEEEmembership{Member,~IEEE}, Wing Shing Wong \IEEEmembership{Fellow,~IEEE} and Huanshui Zhang \IEEEmembership{Senior Member,~IEEE}
\thanks{This work was supported in part by the Research Grants Council of the Hong Kong Special Administrative Region under Project GRF 14630915, National Natural Science Foundation of China under Grant 61803224,
and Schneider Electric, Lenovo Group (China) Limited and the Hong Kong Innovation and Technology Fund (ITS/066/17FP) under the HKUST-MIT Research Alliance Consortium.}
\thanks{C. Tan is with the College of Engineering, QuFu Normal University, Rizhao, Shandong 276800, China and the Department of Information Engineering, The Chinese University of Hong Kong, Shatin, N. T., Hong Kong (\small  e-mail: tancheng1987love@163.com).}
\thanks{W. S. Wong is with the Department of Information Engineering, The Chinese University of Hong Kong, Shatin, N. T., Hong Kong (\small  e-mail: wswong@ie.cuhk.edu.hk).}
\thanks{
H. Zhang is with the School of Control Science and Engineering, Shandong University, Jinan, Shandong 250061, China (\small  e-mail: hszhang@sdu.edu.cn). }
}
\maketitle

\begin{abstract}
In this paper, we mainly investigate an integrated system operating under a software defined network (SDN) protocol.
SDN is a new networking paradigm in which network intelligence is centrally administered and data is communicated via channels that are physically separated from those conveying user data.
Under the SDN architecture, it is feasible to set up multiple flows for transmitting control signals to an actuator with high priority for each individual application.
While each flow may suffer random transmission delay, we focus on the stabilization problem under the joint design of the event-driven strategy in actuator and the control policy in decision-maker.
By introducing a predefined application time, the integrated system can be reformulated as the form of stochastic system with input delay and multiplicative noise.
For such system, we propose a set of necessary and sufficient stabilization conditions.
Specifically, for the scalar system, we derive the allowable sampling period bound that can guarantee stabilization in terms of the probability distributions of the random transmission delays.
A simple example is included to show the performance of our theoretic results.
\end{abstract}


\section{Introduction}

Networked control systems (NCSs) are spatially distributed feedback systems, wherein sensors, controllers and actuators are interconnected through some form of communication networking.
The recent emergence of software defined network (SDN) protocol holds promises for the development of time critical networked control applications among a host of other beneficial attractions; see \cite{wong1997}-\cite{SDN03} for a partial list of references.
SDN is a new type of networking paradigm in which network intelligence is centrally administered and networked control data and decisions are communicated via channels that are physically separated from those conveying user data.
This new architecture dramatically simplifies network management, reduces network latency and opens up access for network control innovation to end users.
Under the SDN architecture, it is feasible to specify individual routing path or even a set of paths for any source and destination node pair at the user application layer.

Due to limited channel capacity, in an open communication network, signal transmission over a routing path invariably experiences random delay or even data packet dropout.
These uncertainties may destabilize the whole system and complicate system analysis.
A wide range of research has been reported dealing with problems related to these uncertainties.
Since it is common to model network-induced delay in terms of probability distribution, some sufficient stabilization conditions were derived for an NCS with random delays; see \cite{nilsson1998}-\cite{fridman2016}.
In \cite{zhanglq2005}, by modelling the sensor-to-controller delay and controller-to-actuator delay as two homogeneous Markov chains, the closed-loop system was reduced to a jump system with two modes while the necessary and sufficient stability conditions were obtained.
On the other hand, NCSs have to deal with the problems of packet dropout.
The possibly simplest model assumes that packet dropouts are the sample-path realization of a Bernoulli process \cite{Sinopoli2004,hu2007} or a homogeneous Markov chain \cite{huang2007,wu2015}.
For example in \cite{elia2005,xiao2012}, the stabilization problem of an NCS over fading channels was considered.
Necessary and sufficient stabilization conditions and an explicit formula for the maximum packet dropout rate were respectively derived.

In this paper, we focus on an integrated model that incorporates both system dynamics and communication network dynamics.
Working under the SDN framework, it is feasible the NCS to establish multiple routing paths for transmitting sampled state information to the decision-maker side and the control signals to the actuator side.
We further assume that each path suffers a random transmission delay that satisfies an identical probability distribution.
Different from most previous studies, our centralized design spaces include the event-driven strategy on the actuator side and the stabilization control policy on the decision-maker side.

In wireless networks, idle-time-based and Markov decision process (MDP)-based scheduling policies were introduced with quality of service (QoS) constraints, where a packet is removed from the system if it is not delivered by the end of the predefined period \cite{kumar,denglei}.
Thus, the scheduling policy can guarantee the delay of each delivered packet is less than one period.
Motivated by the idle-time-based scheduling policy in \cite{kumar}, we present a new event-driven strategy.
By setting a predefined application time, we simplify the sampled system into a discrete-time stochastic system with both input delay and packet dropout.
Based on the Riccati-type and Lyapunov-type approaches, we derive the necessary and sufficient conditions for stabilization.
Motivated by \cite{henrion,chesi}, we also derive a matrix polynomial based condition for stabilizing the uncertain system.
Moreover, for the scalar case, we focus on the allowable sampling period bounds, which are shown to be determined by the system parameter and the probability distribution of random delay.
Note that results in this paper find application in controlling autonomous vehicles over an open access communication network for the purposes of remote site recognition or security control among others.

{\em Notation}:
Let $\mathbb{N}$ denote the set of nonnegative integers, i.e., $\mathbb{N}\triangleq\{0,1,\cdots\}$ and $\mathbb{R}^n$ denote $n$-dimensional real Euclidean space.
Let $\mathbb{S}^{n}$ be the space of all $n$-dimensional symmetric matrices.
For any $q\in \mathbb{R}$, $\lceil q\rceil$ denotes the least integer larger than or equal to $q$ and $\lfloor q\rfloor$ the greatest integer less than or equal to $q$.
$A'$ denotes the transpose of matrix or vector $A$ and $I$ the identity matrix of appropriate dimensions.
$A>0~(\geq 0)$ means that $A$ is a real symmetric positive definite (semidefinite) matrix.
For a matrix polynomial $P(s)$, ${\bf deg}(P(s))$ denotes the maximum of the degree of the entries $P_{i,j}(s)$ in $P(s)$.
$\{w_k,k\in\mathbb{N}\}$ means a sequence of real random variables defined on the
filtered probability space $(\Omega,
\mathcal{F},\mathcal{P};\mathcal{F}_k)$ with the filter $\mathcal{F}_0=\{\emptyset,\Omega\}$
and $\mathcal{F}_k=\sigma\{w_j|j=0,1,2,\cdots,k\}$. Furthermore,
$\hat{x}_{j|k}=E[x_j|\mathcal{F}_k]$ defines the conditional expectation of the state $x_j$
with respect to $\mathcal{F}_k$.

\section{Problem Formulation}

In this paper, we focus on an NCS operating under the SDN protocol as depicted in Fig. 1.
In this model, the plant has a sensor and an actuator, and the decision-maker remotely controls the underlying continuous-time dynamic system:
\begin{align}
\label{sys01}
\dot{x}(t)=Ax(t)+Bu(t),
\end{align}
where $x(t)\in\mathbb{R}^{n_x}$ is the state vector and $u(t)\in\mathbb{R}^{n_u}$ is the
control input vector at time $t\geq 0$.
An open communication network is responsible for the transmission of the sampled state data as well as the control signals.
Without loss of generality, it is assumed that the sensor of the plant is clock driven that periodically samples the state information with a sampling period $h>0$.
The decision-maker is collocated with the sensor and can access the sampled state information at any time slot $\bar{k}=kh$.
Moreover, the actuator of the plant is assumed to be event driven and its event-driven strategy $\eta$ is to be determined.
\begin{figure}[thpb]
  \centering
  \includegraphics[scale=0.6]{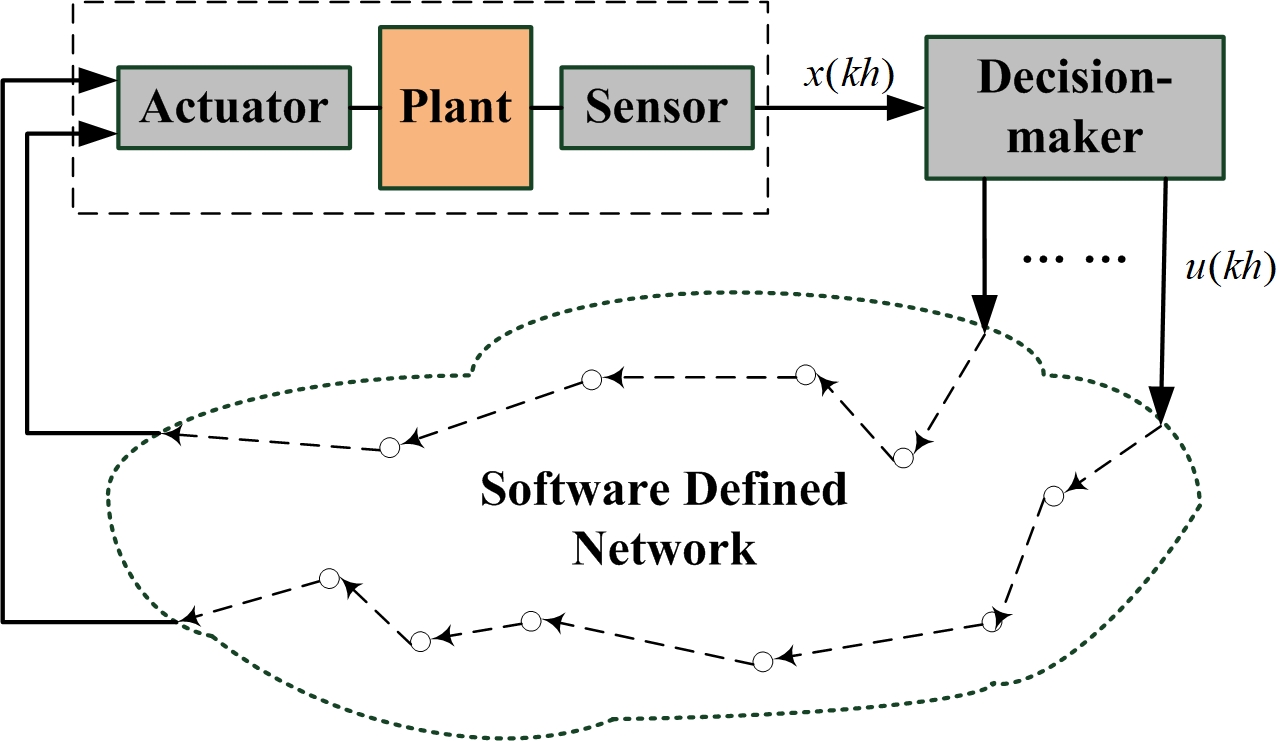}
    \caption{Networked System over SDN}
\end{figure}

Under the SDN architecture, it is feasible to set up multiple information flows for an application operating over a given source destination node pair.
To achieve this, we assume that at the beginning of each time slot $\bar{k}$, the designed control signal with a time stamp is replicated and simultaneously transmitted to the actuator via distinct routing paths, where the number of the multiple flows is given a priori as $m\geq 1$.\
Due to routing and transmission delays, we assume that the transmitted control signal suffers a random transmission delay $d^i_k>0,~k\geq0$, in the $i$-th information flow, $i\in[m]\triangleq\{1,2,\cdots,m\}$.
For the sake of analysis, we further assume that $d_k^i$ is independent of $d_{l}^j$, $k\neq l,~i\neq j$, and $d_k^i$ follows an identical distribution, $F_i$, i.e.,
\begin{align}
F_i(x)=\mathcal{P}(d^i_k\leq x),~i\in[m].
\end{align}

In this paper, we focus on the mean square stabilization problem of system (\ref{sys01}) with a sampling period $h>0$.
To make the underlying system stabilizable, our centralized design spaces include two parts:

$\bullet$ {\em the event-driven strategy $\eta$ on the actuator side};

$\bullet$ {\em the control policy $u_{k}=u(kh)$ on the decision-maker side}, which depends on the proposed event-driven strategy and the current information.

Note that both the event-driven strategy and the control policy influence the dynamics of system (\ref{sys01}).
The problems to be solved are formulated as follows.

{\em Problem I}:
Joint design of the event-driven strategy on the actuator side and the control policy on the decision-maker side to stabilize the sampled system (\ref{sys01}) in the mean square sense.
Explore the necessary and sufficient stabilization conditions.

{\em Problem II}:
For the scalar system and the decoupled multi-input system, derive the allowable sampling period bounds to guarantee stabilization.

\section{Main Results}

\subsection{Event-driven Strategy }

To begin with, we propose a viable event-driven strategy $\eta$ which models the sampled version of the networked control system as a stochastic system.
The basic discriminator of the event-driven strategy is that the control signals need to be delivered to the actuator before a predefined bounded application time $\bar{d}$.

As shown in Fig. 2, we are in a position to describe the event-driven strategy $\eta$ as follows:
By prior agreement, if the sampled control policy $u(kh)$ is not successfully delivered by the end of the predefined application time, it is marked as expired and removed from the SDN.
Otherwise, the actuator of the plant applies the first arriving control signal at a predefined application time $kh+\bar{d}$, where $\bar{d}=dh>0$ is an integral multiple of the sampling period $h>0$ and $d\in\mathbb{N}^+\triangleq\{1,2,\cdots\}$.
That is to say, if the $i$-th transmission delay satisfies $d^i_k\leq \bar{d}$, the first arriving control signal will be applied to the dynamic system at the predefined application time $kh+\bar{d}$; otherwise, all copies of the transmitted control signals are considered as packet dropout.
\begin{figure}[thpb]
  \centering
  \includegraphics[scale=0.6]{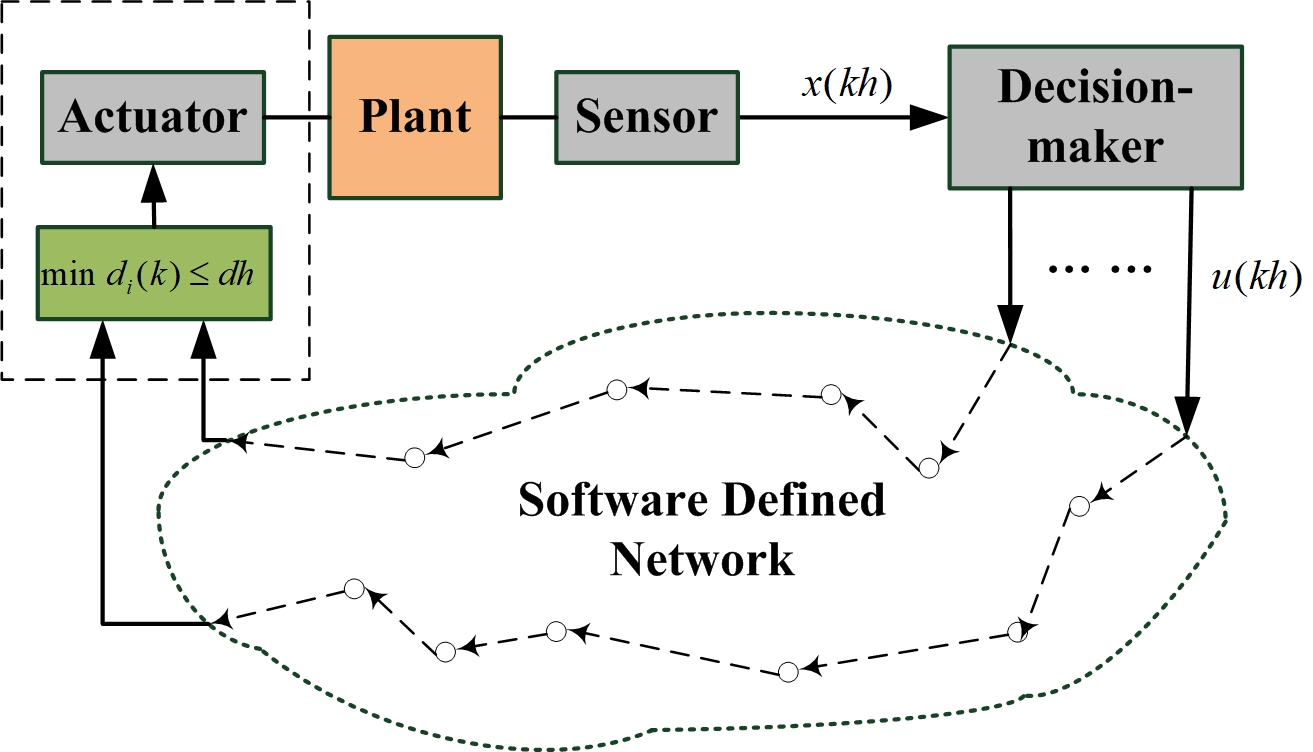}
  \caption{Integrated System with Event-driven Strategy}
\end{figure}

With such event-driven strategy $\eta$, each control signal has a hard delay $dh$ and it becomes useless if it cannot be delivered before the predefined deadline.
By simultaneously considering the random transmission delay and the predefined application time, at any time slot $kh$, the control policy applied by the actuator can be modelled as:
\be
\gamma_{k-d} u((k-d)h),~k\geq 0.
\ee
Here, $\gamma_{k}$ is a binary random variable, where $\gamma_{k}=1$ means that the first arriving data packet has been successfully applied by the actuator and $\gamma_{k}=0$ signifies a dropout of the control signal.
Note that the packet dropout process is caused by the random transmission delay that exceeds the predefined application time.
Based on such simplification, the overall integrated system on the actuator side can be described as:
\be\label{sys02}
x_{k+1}={A}_h x_{k}+\gamma_{k-d}{B}_h u_{k-d},
\ee
where $A_h=e^{Ah}$, $B_h=\int_{0}^h e^{A\tau}Bd\tau$ and $x_{k}=x(kh)$, $u_{k-d}=u((k-d)h)$.
In this case, by utilizing the event-driven strategy $\eta$, we simplify the original model with multiple random delays into a stochastic system with input delay and packet dropout.

Denote $\omega_{k}=\gamma_{k}-E[\gamma_{k}]=\gamma_{k}-(1-p)$.
Then system (\ref{sys02}) can be rewritten as:
\begin{align}
x_{k+1}=A_h x_k+(1-p) B_h  u_{k-d}+w_{k-d} B_h u_{k-d}. \label{sys003}
\end{align}
In this case,
$\{w_k,k\in \mathbb{N}\}$ is a sequence of random variables defined on the
filtered probability space $(\Omega,
\mathcal{F},\mathcal{P};\mathcal{F}_k)$ with $E[w_k]=0$ and $E[w_k w_s]=p(1-p)\delta_{ks}$, where
$\delta_{ks}$ refers to the Kronecker function, i.e., $\delta_{ks}=1$ if $k = s$, and $\delta_{ks}=0$ if $k\neq s$.
In fact, the integrated system (\ref{sys02}) is a special case of the discrete-time stochastic system with multiplicative noise and input delay in \cite{zhanghs2014}, where the stabilization control policy is designed to be the feedback of the conditional expectation of the state.
By extending the definition of mean square stabilization in \cite{zhanghs2014} to the sampled systems, we introduce the following definition.
\begin{definition}\label{def01}
System (\ref{sys02}) is said to be asymptotically mean square stabilizable, if there exists a feedback control policy $u_{k-d}=K\hat{x}_{k|k-d-1}$, such that the closed-loop system is asymptotically mean square stable, i.e., for any initial values, the state $x_k$ satisfies $\lim_{t\rightarrow \infty} {\bf E}\|x_k\|^2=0$.
\end{definition}

\begin{remark}
Recently, some studies have concentrated on the integrated stabilization control with multiple uncertainties.
The most popular methods, including the switched approach \cite{yu2008} and the Lyapunov-Krasivskii functional approach \cite{sun2013}, mainly depend on linear matrix inequalities (LMIs) and only sufficient stabilization conditions were available.
In our previous work \cite{tan2015}-\cite{tan2017}, we focused on an NCS with constant  transmission delay and packet dropout.
We proposed a class of necessary and sufficient stabilization conditions in terms of the positive solution to the delay-dependent algebraic Riccati equation (DARE) or the delay-dependent Lyapunov equation (DLE),  which can be verified as a LMI feasibility problem.
\end{remark}

\subsection{Stabilization Conditions}
For the considered system, based on the DARE and DLE approach in our previous work \cite{tan2017}, we derive a set of stabilization conditions as follows.

\begin{lemma}
\label{lemma01}
Suppose that we apply the event-driven strategy $\eta$ with the given application time $\bar{d}=dh$. The following statements are equivalent.

1) System (\ref{sys02}) is stabilizable in the mean square sense.

2) For any $Q>0$ and $R>0$, there exists a unique positive definite solution $P>0$ satisfying the following DARE:
\begin{equation}
\label{dare01}
P=A_h'PA_h+Q-L'\Phi^{-1}L,
\end{equation}
with
\begin{align}
\label{gare02}
\Phi=&~(1-p)^2 B'_h P B_h+p(1-p)B'_h (A'_h)^d P A^d_h  B_h \nonumber \\
&~+\sum^{d}_{i=0} p(1-p) B'_h (A'_h)^i Q A^i_h  B_h+R ,  \\
L=&~(1-p) B'_h P A_h .
\end{align}
In this case, the stabilization control policy is
\be
u_{k}=-\Phi^{-1}L \hat{x}_{k+d|k-1},
\ee
where
\be
\hat{x}_{k+d|k-1}=A_h^d x_{k} +(1-p) \sum_{i=0}^{d-1} A^i_h B_h u_{k-i-1},
\ee
is the predicted state in time $k+d$ based on the observation of the state at time $k$ and the past control polices.

3) For any $Q>0$, there exist matrices $K$ and $P>0$ satisfying the following DLE:
\begin{align}
P=&~Q+\left[A_h+(1-p)B_hK\right]'P \left[A_h+(1-p)B_hK \right] \nonumber \\
&+p(1-p)K'B_h' (A_h')^d PA_h^dB_hK.
\label{dle}
\end{align}

\end{lemma}

\emph{Proof:} By utilizing the event-driven strategy $\eta$ to the sampled system (\ref{sys01}), the original model is reduced to the integrated system (\ref{sys02}).
Since $\{d_k^i\}_{k\geq 0}$ is assumed to be independent and follow an identical distribution, $\{\gamma_{k}\}_{k\geq 0}$ is an independent and identically distributed Bernoulli process with
\be\label{rate}
\mathcal{P}(\gamma_{k}=0)=p,~\mathcal{P}(\gamma_{k}=1)=1-p
\ee
where $p\in[0,1]$ is the packet dropout rate.
By Theorems 1-2 in \cite{tan2017}, the necessary and sufficient stabilization conditions are derived directly.
\hfill $\blacksquare$

The maximum packet dropout rate is usually introduced to gauge the degree of packet loss that can be tolerated by the integrated system. To guarantee the existence of the maximum packet dropout rate of system (\ref{sys02}), we assume:

{\em
H1) ${A}$ is unstable and its eigenvalues are real and distinct.

H2) ${B}$ has full-column rank.

H3) $({A},{B})$ is a controllable pair.}

\begin{lemma}
\label{lemma02}
Under assumptions H1)-H3), there exists a unique maximum packet dropout rate $p_{\max}\in(0,1]$ such that for any $p<p_{\max}$, system (\ref{sys02}) is stabilizable in the mean square sense.
\end{lemma}

\emph{Proof:}
Under assumptions H1)-H3), for any sample period $h>0$, we obtain that $A_h=e^{Ah}$ is unstable and $B_h=\int_{0}^h e^{A\tau}Bd\tau$ has full-column rank.
Moreover, the discrete-time system $(A_h,B_h)$ is controllable and stabilizable.
Then, by Theorem 1 in \cite{tan2015}, there exists a unique maximum packet dropout rate to guarantee the stabilization.
\hfill $\blacksquare$

Next, we consider a uncertain system where $p$ is time-invariant but uncertain within the interval $[0,\hat{p}]$ and $\hat{p}<p_{\max}$.
For any matrix polynomial $X(p)$, we define parameter-dependent Lyapunov operator $\mathcal {L}_K(\cdot)$ from $\mathbb{R}\times\mathbb{S}^n$ to $\mathbb {S}^n$ as:
\begin{align}
&\mathcal {L}_K(p,X(p)) \nonumber \\
\triangleq &~
X(p)-\left[A_h+(1-p)B_hK\right]'X(p)\left[A_h+(1-p)B_hK\right] \nonumber \\
&-p(1-p)K'B_h'(A_h')^dX(p)A_h^dB_hK.
\label{op01a}
\end{align}
With $\xi=\frac{\hat{p}}{p}-1\geq 0$, we further denote
\begin{align}
X_p(\xi)=&~X(\frac{\hat{p}}{\xi+1}), \label{op01b}\\
\mathcal{H}_K(\xi,X_p(\xi))=&~\mathcal {L}_K(\frac{\hat{p}}{\xi+1},X_p(\xi)). \label{op01c}
\end{align}
Moreover, if a matrix polynomial $X(p)$ can be decoupled to a sum of squares of matrix polynomials, $X(p)$ is said to be a SOS matrix polynomial.
Clearly, if $X(p)$ is a  SOS matrix polynomial, $X(p)$ is positive semidefinite.

Motivated by \cite{henrion,chesi}, we derive the following necessary and sufficient stabilization condition based on matrix polynomials.

\begin{theorem}\label{thm00}
Suppose that we apply the event-driven strategy $\eta$ with the given application time $\bar{d}=dh$.
Under assumptions H1)-H3), the uncertain system (\ref{sys02}) is stabilizable for any $p\in [0,\hat{p}]$ if and only if there exist a matrix polynomial $P(\xi)$ and scalars $\theta>0$, $\zeta\geq  {\bf deg}(P(\xi))$ such that
\begin{align}
P(\xi^2)-(\xi^2+1)^\zeta \theta  I~~~~~~~ &~\text{ is SOS}, \label{sos01} \\
(\xi^2+1)\mathcal{H}_K(\xi^2,P(\xi^2))-(\xi^2+1)^{\zeta+1} \theta I &~\text{ is SOS}. \label{sos02}
\end{align}
\end{theorem}

\emph{Proof:} See Appendix~\ref{app00}. \hfill $\blacksquare$

Note that the LMI feasibility test can be used to verify whether a matrix polynomial is SOS;
see the detail in section \uppercase\expandafter{\romannumeral 3}-A in \cite{chesi}.

\subsection{Sampling Period Bound}
In this part, we study the sampling period bound.
First we focus on the {\bf scalar} case and state the follow assumption:

{\em
H4) $A\geq 0$ and $B\neq0$.
}

In the current scalar model, assumption H4) is equivalent to the two conditions H1) and H2), which also ensure the controllability of system $({A},{B})$ , i.e., assumption H3) holds.
Instead of DARE, DLE and matrix polynomials, we simplify the necessary and sufficient stabilization conditions into a computation formula in terms of one simple inequality.

\begin{theorem}\label{thm01}
Under assumption H4), system (\ref{sys02}) is stabilizable in the mean square sense if and only if
there exists an application time $\bar{d}=dh>0$ such that
\begin{equation}
\prod_{i=1}^m (1-F_i(\bar{d}))<\frac{1}{{e}^{2A\bar{d}}(e^{2Ah}-1)+1}.\label{inequ01}
\end{equation}
\end{theorem}

\emph{Proof:}
According to the above definition of the event-driven strategy $\eta$, the occurrence of packet dropout means that all copies of the control signals are not successfully delivered to the actuator by the end of the predefined application time.
Based on the assumption of $\{d_k^i\}_{k\geq 0}$, it follows that the packet dropout rate  $p\in[0,1]$ in (\ref{rate}) satisfies
\be
p=\mathcal{P}(\min_{1\leq i\leq m}d^i_k>\bar{d})=\prod_{i=1}^m (1-F_i(\bar{d})).\label{packet}
\ee
In \cite{tan2015}, the maximum packet dropout rate of the integrated system (\ref{sys02}) is derived as
\be
p_{\max}=\frac{1}{e^{2Ah(d+1)}-e^{2Ahd}+1}.
\ee
System (\ref{sys02}) is stabilizable in the mean square sense if and only if $p<p_{\max}$, i.e., the inequality (\ref{inequ01}) holds.
\hfill $\blacksquare$

In date network, the random transmission delay depends on the service time of the single server \cite{rap1996}.
Consider the special case that the delays in different paths can be modelled by independent exponential distributions.
We put forward the following assumption:

{\em
H5)
For any $i\in[m]$, the random transmission delay $d_k^i$ satisfies the following probability distribution
\be
F_i(x)=1-e^{-r_i x},
\ee
where $r_i>0$ reflects service reliability of the $i$-th information flow.}

Define $\bar{r}=\sum_{i=1}^m r_i$ as the total service rate of all routing paths, which influences the stabilization of system (\ref{sys02}).
Next we study the allowable sampling period bound of the integrated system (\ref{sys02}).

\begin{theorem}\label{thm02}
Under assumptions H4)-H5), we have the following statements.

1) Suppose $\bar{r}>2A$. For any bounded sampling period $0<h<\infty$, there exists a positive integer $d>0$ such that system (\ref{sys02}) is stabilizable in the mean square sense.

2) Suppose $\bar{r}=2A$. For any sampling period $0<h<\bar{h}$, there exists a positive integer $d>0$ such that system (\ref{sys02}) is stabilizable in the mean square sense, where $\bar{h}$ is the explicit sampling period bound with
\begin{align}
\bar{h}=\frac{\ln2}{2A}.
\end{align}

3) Suppose $0<\bar{r}<2A$.
For any bounded sampling period $h_u \leq h<\infty$, system (\ref{sys02}) cannot be stabilized  for any $d>0$, where $h_{u}$ is the upper sampling period bound satisfying
\begin{align}
h_{u}=\frac{1}{2A}\ln[\frac{\bar{r}(2A-\bar{r})^{\frac{2A}{\bar{r}}-1}}{(2A)^{\frac{2A}{\bar{r}}}}+1].
\end{align}
For any sampling period $0<h\leq h_{l}$, there exists a positive integer $d>0$ such that system (\ref{sys02}) is stabilizable in the mean square sense, where $h_{l}>0$ is the lower sampling period bound satisfying
\begin{align}
h_{l}=&~\max\{h_{l_1},h_{l_2}\}, \\
h_{l_{1}}=& ~\frac{\ln\left(
e^{(\bar{r}-2A)(\bar{t}+ h_u)}-e^{-2A(\bar{t}+ h_u)}+1\right)}{2A}, \\
h_{l_{2}}=& ~\frac{\ln\left(
e^{(\bar{r}-2A)(\bar{t}- h_u)}-e^{-2A(\bar{t}- h_u)}+1\right)}{2A}, \\
\bar{t}=&~\frac{\ln(2A)-\ln(2A-\bar{r})}{\bar{r}}. \label{def-t}
\end{align}
\end{theorem}

\emph{Proof:} See Appendix~\ref{app01}. \hfill $\blacksquare$

In the proof of Theorem \ref{thm02}, we obtain $h_{l}<h_u$ from (\ref{re01}) and (\ref{re02}).
When $0<\bar{r}<2A$, we classify the sampling period according to two sampling period bounds.
For $h_l<h\leq h_u$, we derive a stabilization criterion as follows.

\begin{proposition}\label{proposition01}
Suppose $0<\bar{r}<2A$.
Under assumptions H4)-H5), for any sampling period $h_l<h\leq h_u$, if the following inequality holds
\begin{align}
e^{2Ah}-1<f^*(d),
\end{align}
where
\begin{align}
f^*(d)=&~\max\{f(d_1),f(d_2)\},\\
f(d_i)=&~(e^{\bar{r}d_ih}-1)e^{-2Ad_ih},~i=1,2, \\
d_1=&~ \lceil \frac{\bar{t}}{h}\rceil,~d_2= \lfloor \frac{\bar{t}}{h}\rfloor,
\end{align}
then system (\ref{sys02}) is stabilizable in the mean square sense. 
\end{proposition}

\emph{Proof:}
For any sampling period $h>0$ satisfying $h_l<h\leq h_u$, denote $d_1= \lceil \frac{\bar{t}}{h}\rceil$ and $d_2= \lfloor \frac{\bar{t}}{h}\rfloor$. 
It follows that
\begin{align}
\bar{t}-h<d_2 h\leq \bar{t}\leq d_1 h<\bar{t}+h.
\end{align}
The maximum point for $f(d)=e^{(\bar{r}-2A)dh}-e^{-2Adh}$ is
\begin{align}
f^*(d)=\max_{i=1,2}\{(e^{\bar{r}d_ih}-1)e^{-2Ad_ih} \}.
\end{align}
By (\ref{inequ01}), if
\be
e^{2Ah}-1<f^*(d),
\ee
system (\ref{sys02}) is stabilizable in the mean square sense with the application time $\bar{d}_{\bar{i}}=hd_{\bar{i}}$, where $\bar{i}\in\{1,2\}$ and $f(d_{\bar{i}})=f^*(d)$.
\hfill $\blacksquare$

In principle, by means of Theorem \ref{thm02} and Proposition \ref{proposition01}, we present the stabilization criteria for any sampling period $h>0$.
It is shown that the stabilization region (the sampling period bounds) depends on the unstable parameter $A\geq 0$ and the service capability $\bar{r}>0$.
When $\bar{r}>0$ is given a priori, the stabilization region grows as $A\geq 0$ decreases to zero.
Suppose system (\ref{sys01}) is given a priori, i.e., $A$ is fixed.
Note that the stabilization region grows with increasing service capability $\bar{r}>0$.
These theoretic results can be illustrated in Table I.
\begin{table}[thpb]
\centering
\caption{Sampling Period Bounds with $A$ and $\bar{r}$}
  \begin{tabular}{|c|c|c|c|c|c|c|c|}
  \hline
   $A$ & $\bar{r}$ & $h_u$ & $h_l$ &  $A$ & $\bar{r}$ & $h_u$ & $h_l$\\
       \hline
   $0.6$ & $1$ & $0.3824$ & $0.3779$  & $0.4$ & $0.1$ & $0.0599$ & $0.0598$ \\
       \hline
   $0.7$ & $1$ & $0.2569$ & $0.2534$ & $0.4$ & $0.2$ & $0.1253$ & $0.1249$ \\
    \hline
   $0.8$ & $1$ & $0.1862$&  $0.1838$ & $0.4$ & $0.3$ & $0.1977$ & $0.1964$ \\
    \hline
   $0.9$ & $1$ & $0.1416$ &  $0.1401$ & $0.4$ & $0.4$ & $0.2789$ & $0.2761$ \\
       \hline
   $1$ & $1$ & $0.1116$ &  $0.1105$   & $0.4$ & $0.5$ & $0.3723$ & $0.3676$ \\
       \hline
   $1.1$ & $1$ & $0.0902$ &  $0.0895$ & $0.4$ & $0.6$ & $0.4837$ & $0.4773$ \\
       \hline
   $1.2$ & $1$ & $0.0745$ &  $0.0740$ & $0.4$ & $0.7$ & $0.6261$ & $0.6169$ \\
       \hline
\end{tabular}
\end{table}

To explore the nature of the vector system, we focus on a special decoupled system (\ref{sys02}) with
\bee\label{eq201}
A=diag\{A_1,A_2,\cdots,A_n\},~B=diag\{B_1,B_2,\cdots,B_n\},
\eee
where $A_i,B_i\in\mathbb{R}$ are scalars.
In this case, assumptions H1)-H2) become that $A$ is unstable and $rank(B)=n$,
i.e., there exists at least a integer $i\in [n]$ such that $A_i\geq 0$ and for any $i\in[n]$, $B_i\neq 0$.
For convenience, we further assume that
\be\label{inequ02}
 A_1\geq \cdots \geq A_\mu\geq 0>A_{\mu+1}\geq \cdots \geq A_{n},~\mu\leq n.
\ee
In this case,  $({A},{B})$ is a controllable pair.
Next we are in a position to present the following stabilization result. 

\begin{theorem}\label{thm03}
Under assumption H1)-H2), system (\ref{sys02}) is stabilizable in the mean square sense if and only if there exists an application time $\bar{d}=dh>0$ such that
\begin{equation}\label{inequ03}
\prod_{i=1}^m (1-F_i(\bar{d}))<\frac{1}{{e}^{2A_1\bar{d}}(e^{2A_1h}-1)+1}.
\end{equation}
\end{theorem}

\emph{Proof:}
It follows from Theorem \ref{thm01} that the $i$-th scalar system (\ref{sys02}) is stabilizable with $A_i,B_i$ if and only if the packet dropout rate $p\in[0,1]$ in (\ref{packet}) satisfies
\begin{equation}
p<\frac{1}{{e}^{2A_i\bar{d}}(e^{2A_ih}-1)+1}.
\end{equation}
It follows that the decoupled system (\ref{sys02}) is stabilizable if and  only if
\begin{equation}
p<\min_{i\in[\mu]}\frac{1}{{e}^{2A_i\bar{d}}(e^{2A_ih}-1)+1},
\end{equation}
which implies (\ref{inequ03}) based on (\ref{inequ02}) and completes the proof.
\hfill $\blacksquare$

\begin{remark}
Note that the stabilization of the decouple system (\ref{sys02}) only depends on the maximum unstable eigenvalue $\lambda_i=A_i\geq 0$.
Similar to the above discussion, we can derive the sampling period with $A_1$.
For the general vector system, it is an open question on how to calculate the sampling period bounds.
\end{remark}

\begin{remark}
In this paper, we assume that there is no traffic from the sensor to the decision-maker.
From a general application perspective, on the decision-maker side, it is common to model the arrival state information as a random process with intermittent observations; see \cite{Sinopoli2004,qingyuan}.
The stabilization problem of such system under the SDN protocol defines a promising and challenging research direction.
\end{remark}
\section{Simulation}
In this section, we give a simple example to illustrate our main results.
Consider the scalar system  with $A=0.25,~B=1,~d=2,~m=2,$ and the initial condition $x_0=2$, $u_{-2}=u_{-1}=0$.

If $r_1=r_2=0.5$, then $\bar{r}=r_1+r_2>2A$. By Theorem \ref{thm02}, for any bounded sampling period $h>0$, the integrated system (\ref{sys02}) is stabilizable for some application time $\bar{d}$.
In this case, when $h=1$, we have $\bar{A}=e^{Ah}=1.2840$, $\bar{B}=\int_0^h e^{As}Bds=1.1360$ and $p=e^{-r_1mdh}=0.1353$.
Solving the DARE (\ref{dare01}) implies that the unique positive solution is $P=3.7959$.
By Lemma \ref{lemma01}, the stabilization control policy can be designed as
$u^1_k=K_P \hat{x}_{k|k-3}=-0.7232 \hat{x}_{k|k-3}.$
The simulation of the state response is shown with the red line in Fig. 3.

If $r_1=r_2=0.2$, then $\bar{r}=r_1+r_2<2A$. By Theorem \ref{thm02}, we have
\begin{align}
h_{u}
&=\frac{1}{2A}\ln[\frac{\bar{r}(2A-\bar{r})^{\frac{2A}{\bar{r}}-1}}{(2A)^{\frac{2A}{\bar{r}}}}+1]=0.8571.
\end{align}
In this case, for sampling period $h=1>h_{u}$, system (\ref{sys02}) cannot be stabilized with any controller. The simulation of the state response $E(x^2_k)$ with $u^1_k=K_P \hat{x}_{k|k-3}$ is shown with the blue line in Fig. 3.

To study the conservation of the developed sampling period bounds, we define $\Delta h=h_u-h_l$.
For different service capability $\bar{r}$, the sampling period bounds can be summarized in Table I, which shows the validity of our theoretic results.
\begin{table}[thpb]
\centering
\caption{Sampling Period Bounds}
  \begin{tabular}{|c|c|c|c|c|}
  \hline
   $A$ & $\bar{r}$ & $h_u$ & $h_l$ & $\Delta h=h_u-h_l$\\
       \hline
   $0.5$ & $0.9$ & $0.5288$ & $0.5240$ & $0.0048$  \\
       \hline
   $0.5$ & $0.75$ & $0.3869$ & $0.3828$ & $0.0041$  \\
       \hline
   $0.5$ & $0.6$ & $0.2820$ & $0.2785$ & $0.0035$  \\
    \hline
   $0.5$ & $0.45$ & $0.1962$ & $0.1944$ & $0.0018$  \\
    \hline
   $0.5$ & $0.3$ & $0.1227$ & $0.1221$ & $0.0006$  \\
       \hline
   $0.5$ & $0.15$ & $0.0580$ & $0.0579$ & $0.0001$  \\
\hline
\end{tabular}
\end{table}

\section{Conclusion}

In this paper, we focus on the stabilization problem of the integrated system over SDN. By utilizing a new event-driven strategy on the actuator side, we present a set of the necessary and sufficient stabilization conditions are present with DARE, DLE and matrix polynomials.
Moreover, for the scalar case, we derive the allowable sampling period bounds.

In the future, it is promising to generalize our system in several aspects.
On the one hand, we have only proposed the sampling bounds for the scalar case and thus the general vector case is worth considering.
On the other hand, our results are based on the event-driven strategy $\eta$, which simplifies the integrated system into a stochastic system with constant input delay and packet dropout but shrinks the design space of control policy.
It would be challenging to figure out more event-driven strategies, control polices and analysis framework for sampling period bounds.

\appendices
\section{Proof of Theorem \ref{thm00}}\label{app00}
\textbf{Proof.}
\emph{Sufficiency}.
It follows from Theorem 4 in \cite{chesi} that (\ref{sos01})-(\ref{sos02}) hold if and only if the following inequalities hold for any $\xi\in\mathbb{R}$,
\begin{align}
&~~~~~~~~~~~~~P(\xi^2)-(\xi^2+1)^\zeta  \theta I \geq 0, \label{sos0a}\\
&(\xi^2+1)\mathcal{H}_K(\xi^2,P(\xi^2))-(\xi^2+1)^{\zeta+1}\theta I \geq 0, \label{sos0b}
\end{align}
which implies that
\begin{align}
P(\xi)>0,~\mathcal{H}_K(\xi,P(\xi))>0,~\forall \xi\geq 0. \label{su01}
\end{align}
By applying $\xi=\frac{\hat{p}}{p}-1\geq 0$ to (\ref{su01}), we obtain that for any $p=\frac{\hat{p}}{\xi+1}\in [0,\hat{p}]$, there exists a positive definite matrix $X(p)=P(\frac{\hat{p}}{p}-1)>0$ satisfying $\mathcal {L}_K(p,X(p))>0$.
By Theorem 2 in \cite{tan2017}, we have that system (\ref{sys02}) is stabilizable in the mean square sense.

\emph{Necessity}.
Suppose system (\ref{sys02}) is stabilizable for any $p\in[0,\hat{p}]$.
Based on Theorem 2 in \cite{tan2017}, the following delay-free stochastic system is stabilizable in the mean square
\begin{align}
z(k+1)=&~A_h z(k)+(1-p)B_h v(k)+w_k A^d_hB_h v(k), \label{sys3b}
\end{align}
where $v(k)=Kz(k)$ is a stabilization control policy and $w_k$ follows the same probability distribution in (\ref{sys003}).
It follows from Theorem 1 in \cite{zhouxy} that there exists a unique positive definite matrix $X(p)>0$ that satisfies the DLE $\mathcal {L}_K(p,X(p))=I$ in (\ref{dle}).
Rewrite $\mathcal {L}_K(p,X(p))$ as follows
\begin{align}
&\mathcal {L}_K(p,X(p)) \nonumber \\
=&~
X(p)-(1-p)\left[A_h+B_hK\right]'X(p)\left[A_h+B_hK\right] \nonumber \\
&-pA_h'X(p)A_h+p(1-p)K'B_h'X(p)B_hK  \nonumber \\
&-p(1-p)K'B_h'(A_h')^dX(p)A_h^dB_hK.
\end{align}
which indicates that $\mathcal {L}_K(p,X(p))$ is linear with respect to $X(p)$
and $\mathcal {L}_K(p,X)$ is a matrix polynomial in $p\in[0,\hat{p}]$. It follows from \cite{su2017} that $X(p)$ is a rational function of $p$ which can be rewritten as
$X(p)=Y(p)/y(p)$, where $Y(p)>0$ is a matrix polynomial and $y(p)>0$ is a polynomial.
There exists a scalar $\theta>0$ such that
\begin{align}
Y(p)-\theta I\geq 0, ~\mathcal {L}_K(p,Y(p))-\theta I\geq 0,~\forall p\in[0,\hat{p}].
\label{sos03}
\end{align}
Denote $P_p(\xi)=Y(\frac{\hat{p}}{\xi+1})$ with $\xi=\frac{\hat{p}}{p}-1\geq 0$.
It follows that
\begin{align}
(\xi+1)^\zeta P_p(\xi)- (\xi+1)^\zeta \theta I \geq  & ~0, \label{sos05} \\
(\xi+1)^{\zeta+1} \mathcal {H}_K(\xi,P_p(\xi))-(\xi+1)^{\zeta+1}\theta I \geq & ~0,
\forall \xi\geq 0,\label{sos06}
\end{align}
where $\zeta={\bf deg}(Y(p))$ and $\mathcal {H}_K(\cdot)$ is defined in (\ref{op01c}).
In this case, when we denote $P(\xi)=(\xi+1)^\zeta P_p(\xi)$, we obtain
\be
(\xi+1)^{\zeta+1} \mathcal {H}_K(\xi,P_p(\xi))=(\xi+1) \mathcal {H}_K(\xi,P(\xi)).
\ee
Moreover, $P(\xi)$ and $(\xi+1) \mathcal {H}_K(\xi,P(\xi))$ are both matrix polynomials with respect to $\xi$ and ${\bf deg}(P(\xi))\leq \zeta$.
Therefore, (\ref{sos0a})-(\ref{sos0b}) hold for any $\xi\in\mathbb{R}$, which implies that (\ref{sos01})-(\ref{sos02}) hold.
The proof is completed.
\hfill $\blacksquare$

\section{Proof of Theorem \ref{thm02}}\label{app01}
\textbf{Proof.}
Under assumption H5), the packet dropout rate can be formulated as:
\be
p=\prod_{i=1}^m (1-F_i(\bar{d}))=e^{-\bar{r}\bar{d}},~\bar{r}=\sum_{i=1}^m r_i.
\ee
By Theorem \ref{thm01}, by utilizing the event-driven strategy $\eta$, the integrated system (\ref{sys02}) is stabilizable if and only if there exist a sampling period $h>0$ and an application time $d>0$ satisfying
\begin{equation}\label{pro01}
e^{-\bar{r}\bar{d}}<\frac{1}{{e}^{2A\bar{d}}(e^{2Ah}-1)+1}.
\end{equation}
Since ${e}^{2A\bar{d}}(e^{2Ah}-1)>1$ holds for any $h>0$, (\ref{pro01}) is equivalent to
\begin{equation}\label{pro02}
e^{2Ah}-1<e^{(\bar{r}-2A)\bar{d}}-e^{-2A\bar{d}}.
\end{equation}
Define $g(h)=e^{2Ah}-1$. It is evident that $g(h_1)\geq g(h_2)\geq0$ holds for any $h_1\geq h_2>0$.
For any $t>0$, define
\be f(t)=e^{(\bar{r}-2A)t}-e^{-2At}\geq 0. \ee
It follows that
\begin{equation}\label{pro03}
f'(t)=(\bar{r}-2A)e^{(\bar{r}-2A)t}+2Ae^{-2At}.
\end{equation}
For $t_0=0$, we have $f(0)=0.$

1) If $\bar{r}>2A$, $f'(t)>0$ holds for any $t>0$. In this case, we have $\lim_{t\rightarrow\infty}f(t)=\infty$.
For any bounded sampling period $0<h_1<\infty$, there exists a sufficiently large integer $d_1>0$ such that
\be
e^{2Ah_1}-1<f(d_1h_1),
\ee
which indicates that system (\ref{sys02}) is stabilizable with $\bar{d}_1=h_1d_1$.

2) If $\bar{r}=2A$, $f(t)=1-e^{-2At}$ is monotone increasing with $\lim_{t\rightarrow\infty}f(t)=1$.
From (\ref{pro02}), we have
\be
0<h<\frac{\ln2}{2A}=\bar{h}.
\ee
For any sampling period $0<h_2<\bar{h}$, there exists a sufficiently large integer $d_2>0$ such that
\be
e^{2Ah_2}-1<f(d_2 h_2)<1,
\ee
which implies that system (\ref{sys02}) is stabilizable with $\bar{d}_2=h_2d_2$.

3) If $0<\bar{r}<2A$, then $\lim_{t\rightarrow\infty}f(t)=f(0)=0$. Solving $f'(\bar{t})=0$ implies
\begin{equation}
\bar{t}=\frac{\ln(2A)-\ln(2A-\bar{r})}{\bar{r}}>0.
\end{equation}
It is easy to check $f''(\bar{t})<0$, which indicates that $\bar{t}>0$ is the maximum point of $f(t)$ with
\begin{align}
f(\bar{t})
&=\frac{\bar{r}(2A-\bar{r})^{\frac{2A}{\bar{r}}-1}}{(2A)^{\frac{2A}{\bar{r}}}}.
\end{align}
To guarantee the inequality (\ref{pro02}) holds, we have
\begin{equation}\label{pro04}
\begin{aligned}
0<h&<\frac{\ln(f(\bar{t})+1)}{2A}
\\&=\frac{1}{2A}\ln[\frac{\bar{r}(2A-\bar{r})^{\frac{2A}{\bar{r}}-1}}{(2A)^{\frac{2A}{\bar{r}}}}+1]
=h_u.
\end{aligned}\end{equation}
Hence, for any bounded sampling period $h_u<h_3<\infty$ and $d>0$, we have
\begin{equation}
e^{2Ah_3}-1\geq f(\bar{t}) \geq e^{\bar{r}dh_3-2Adh_3}-e^{-2Adh_3},
\end{equation}
which implies that system (\ref{sys02}) cannot be stabilized for any positive  $d>0$.

Denote $\bar{l}$ satisfying
\begin{equation}
\bar{l} = \frac{\bar{t}}{h_u}=\frac{2A\bar{t}}{\ln[e^{(\bar{r}-2A)\bar{t}}-e^{-2A\bar{t}}+1]}.
\end{equation}
It follows that
\begin{align*}
f((\bar{l}+ 1)h_u)=f(\bar{t}+ h_u)<f(\bar{t}), \\
f((\bar{l}- 1)h_u)=f(\bar{t}- h_u)<f(\bar{t}).
\end{align*}
Define
\begin{align}
h_{l_{1}}=& \frac{\ln(f(\bar{t}+ h_u)+1)}{2A} \nonumber \\
=& ~\frac{\ln\left(
e^{(\bar{r}-2A)(\bar{t}+ h_u)}-e^{-2A(\bar{t}+ h_u)}+1\right)}{2A},
\end{align}
which implies that
\be h_{l_1}<h_u,~e^{2A h_{l_1} }-1=f(\bar{t}+ h_u). \label{re01} \ee
In this case, for any sampling period $0<h_4\leq h_{l_{1}}$, define
$d_4= \lceil \frac{\bar{t}}{h_4}\rceil>0$.
Then, we have
\be
\bar{t}\leq \bar{d}_4=d_4h_4<\bar{t}+h_4<\bar{t}+h_u.
\ee
Since $f(t)$ is monotone decreasing for $t\geq \bar{t}$, it follows that
\begin{align}
e^{2Ah_4}-1\leq e^{2Ah_{l_1}}-1=f(\bar{t}+h_u)<f(\bar{d}_4),
\end{align}
which indicates that system (\ref{sys02}) is stabilizable with $\bar{d}_4>0$.
Similarly, define
\begin{align}
h_{l_{2}}
=& ~\frac{\ln\left(
e^{(\bar{r}-2A)(\bar{t}-h_u)}-e^{-2A(\bar{t}- h_u)}+1\right)}{2A},
\end{align}
which implies that
\be h_{l_2}<h_u,~,e^{2A h_{l_2} }-1=f(\bar{t}- h_u). \label{re02}\ee
For any sampling period $0<h_5\leq h_{l_{2}}$, denote $d_5= \lfloor \frac{\bar{t}}{h_5}\rfloor$.
Then, we have
\be
\bar{t}-h_u<\bar{t}-h_5 <\bar{d}_5=d_5 h_5\leq \bar{t}.
\ee
Since $f(t)$ is monotone increasing for $0<t\leq \bar{t}$,
it follows that
\begin{align}
e^{2Ah_5}-1\leq e^{2Ah_{l_2}}-1=f(\bar{t}-h_u)<f(\bar{d}_5),
\end{align}
which indicates that system (\ref{sys02}) is stabilizable with $\bar{d}_5>0$.
Hence, for any sampling period $0<h \leq h_{l}$, there exists some $d>0$ such that system (\ref{sys02}) is stabilizable in the mean square sense and this proof is completed.
\hfill $\blacksquare$

\end{document}